\documentclass[12pt]{amsart}
\usepackage{amscd,amsmath,amssymb,amsfonts}
\theoremstyle{plain}
\newtheorem{thm}{Theorem}
\newtheorem{lem}[thm]{Lemma}

\newtheorem{ex}[thm]{Example}

\numberwithin{thm}{section}
\numberwithin{equation}{section}
\newcommand{\p}{\partial}

\newcommand{\eq}[2]{\begin{equation}\label{#1}#2 \end{equation}}
\newcommand{\ml}[2]{\begin{multline}\label{#1}#2 \end{multline}}

\newcommand{\sB}{{\mathcal B}}
\newcommand{\sC}{{\mathcal C}}
\newcommand{\sD}{{\mathcal D}}
\newcommand{\sE}{{\mathcal E}}

\newcommand{\sO}{{\mathcal O}}

\newcommand{\C}{{\mathbb C}}

\newcommand{\N}{{\mathbb N}}
\renewcommand{\P}{{\mathbb P}}

\newcommand{\Z}{{\mathbb Z}}
\begin{document}

\title[Homology]{Homology for irregular connections} 
\author{Spencer Bloch}
\address{Dept. of Mathematics,
University of Chicago,
Chicago, IL 60637,
USA}
\email{bloch@math.uchicago.edu}

\author{H\'el\`ene Esnault}
\address{Mathematik,
Universit\"at Essen, FB6, Mathematik, 45117 Essen, Germany}
\email{esnault@uni-essen.de}
\date{May 12, 2000}
\begin{abstract}Homology with values in a connection with possibly
  irregular singular points on an algebraic curve is defined,
  generalizing homology with values in the underlying local system for
  a connection with regular singular points. Integration defines a
  perfect pairing between de Rham cohomology with values in the
  connection and homology with values in the dual connection. 
\end{abstract}
\subjclass{Primary 14C40 19E20 14C99}
\maketitle
\begin{quote}

\end{quote}

\section{Introduction}

Consider the following formulas, culled, one may imagine, from a textbook
on calculus:
\begin{xalignat*}{3} && \sqrt{\pi} = \int_{-\infty}^{\infty}e^{-t^2}dt &&
\\
&&(e^{2\pi is}-1)\Gamma(s) = (e^{2\pi is}-1)\int_0^\infty
e^{-t}t^s\frac{dt}{t} && \qquad
\text{     Gamma function} \\
&& J_n(z) = \frac{1}{2\pi i}\int_{\{|u|=\epsilon\}}
\exp(\frac{z}{2}(u-\frac{1}{u}))\frac{du}{u^{n+1}} && \text{ Bessel
function.}
\end{xalignat*}

These are a few familiar examples of periods associated to 
connections with irregular singular points on Riemann surfaces.
Curiously, though of course such integrals have been studied for 200
years or so, and mathematicians in recent years have developed a powerful
duality theory for holonomic $\sD$-modules (for dimension $1$, which is
the only case we will consider, cf.
\cite{M1},chap. IV, and \cite{M2}), it is not easy from the literature
to interpret such integrals as periods arising from
a duality between homological cycles and differential forms. A
homological duality of this sort is well understood for differential
equations with regular singular points. Our purpose in this note is to
develop a similar theory in the irregular case. Of course, most of the
``heavy lifting'' was done by Malgrange op. cit. We hope, in
reinterpreting his theory, to better understand relations between
irregular connections and wildly ramified $\ell$-adic sheaves. There are
striking relations between
$\epsilon$-factors for $\ell$-adic sheaves on curves over finite fields
and determinants of irregular periods \cite{T} which merit further study.
Finally, relations between irregular connections and the arithmetic
theory of motives remain mysterious. 

Let $X$ be a smooth, compact, connected algebraic curve (Riemann surface)
over $\C$. Let $D=\{x_1,\dotsc,x_n\}\subset X$ be a non-empty, finite set
of points (which we also think of as a reduced effective divisor), and
write $U := X\setminus D\stackrel{j}{\hookrightarrow} X$. Let $E$ be a vector
bundle on $X$, and suppose given a connection with meromorphic poles on
$D$
$$\nabla : E \to E\otimes \omega(*D).
$$
Here $\omega$ is the sheaf of holomorphic $1$-forms on $X$ and $*D$
refers to meromorphic poles on $D$. Unless otherwise indicated, we work
throughout in the analytic topology. The de Rham cohomology
$H^*_{DR}(X\setminus D;
E,\nabla)$ is the cohomology of the complex of sections
\eq{0.1}{\Gamma(X, E(*D)) \stackrel{\nabla}{\longrightarrow}
\Gamma(X,E\otimes
\omega(*D))
}
placed in degrees $0$ and $1$. These cohomology groups are finite
dimensional \cite{DeI}, Proposition 6.20, (i). 

Let $E^\vee$ be the dual bundle, and let $\nabla^\vee$ be the dual
connection, so
\eq{0.2}{d\langle e,f^\vee \rangle = \langle \nabla(e),f^\vee \rangle +
\langle e,\nabla^\vee(f^\vee)\rangle 
}
Define $\sE = \ker(\nabla)$, and $\sE^\vee = \ker(\nabla^\vee)$ to be the
corresponding local systems of flat sections on $U$. We want to define
homology with values in these local systems, or more precisely with
values in associated cosheaves on $X$.  For
$x\in X\setminus D$, $\sE_x$ will denote the stalk of $\sE$ at $x$.
Define the co-stalk at $0\in D$
\eq{0.3}{\sE_0 := \sE_x/(1-\sigma)\sE_x
}
where $x\neq 0$ is a nearby point, and $\sigma$ is the local monodromy
about $0$. We write
$\sC_n=\sC_n(E,\nabla)$ for the group of
$n$-chains with values in
$\sE$ and rapid decay near $0$. Write $\Delta^n$ for the $n$-simplex
and $b\in \Delta^n$ for its barycenter. Thus,
$\sC_n(E,\nabla)$ is spanned by elements $c\otimes\epsilon$ with
$c:\Delta^n \to X$ and
$\epsilon
\in \sE_{c(b)}$, where $b\in \Delta^n$ is the barycenter. We
assume 
$c^{-1}(0) =
\text{union of faces}\subset \Delta^n$ and that $\epsilon$ has
{\it rapid decay}
near $D$. 
This is no condition if $D\cap c(\Delta^n)=\emptyset$. If $0\in
D\cap c(\Delta^n)$, we take $e_i$ a basis for $E$ near $0$ and write
$\epsilon = \sum f_i c^*(e_i)$. Let $z$ be a local parameter at $0$ on
$\Delta$. We require that for all $N\in \N$, constants
$C_N>0$ exist with 
$|f_i(z)| \le C_N |z|^{N}$ on $\Delta^n\setminus c^{-1}(0)$. 
Note that if $\nabla$ has logarithmic poles in one point, then
rapid decay implies vanishing. Thus in this case, we deal with
the sheaf $j_!\sE$, where $j: X\setminus D\to X$.

There is a natural boundary map
\eq{0.4}{ \partial: \sC_n(E,\nabla) \to \sC_{n-1}(E,\nabla);\quad
\partial(c\otimes\epsilon) = \sum (-1)^j c_j\otimes\epsilon_j 
}
where $c_j$ are the faces of $c$. Note if $b_j$ is the barycenter of the
$j$-th face and $c(b_j)\neq 0$, $c$ determines a path from $c(b)$ to
$c(b_j)$ which is canonical upto homotopy on $\Delta \setminus \{0\}$. (As a
representative, one can take $c[b_j,b]$, the image of the straight line
from $b$ to $b_j$. By assumption, $c^{-1}(0)$ is a union of faces, so it
does not meet the line.) Thus $\epsilon \in \sE_{c(b)}$ determines
$\epsilon_j \in \sE_{c(b_j)}$. Similarly for $0\in D$, if $c(b_j)=0$ there
is corresponding to $\epsilon$ a unique
$\epsilon_j\in \sE_0$ because we have taken coinvariants. 
If $c:\Delta^n \to D$ is a constant simplex, there is no rapid decay
condition.

It is straightforward to compute that
$\partial\circ\partial = 0$. Consider $c\otimes \epsilon$.  
If $c(b)=0$, where $b\in \Delta^2$ is the barycentre, then
$c(\Delta^2)=0$ and 
$\epsilon=\epsilon_i=(\epsilon_i)_j \in \sE_0$ for all $i$ and $j$
involved, thus the condition is trivially fulfilled. 
If not, and some $c(b_i)=0$, then $(\epsilon_j)_i=(\epsilon_i)_j
\in \sE_0$ for all $j$, and if all $c(b_i)\ne 0$, then one has
by unique analytic continuation in $c(\Delta^2)$ the relation $
(\epsilon_i)_j=(\epsilon_j)_i \in \sE_{{\rm edge}_{ij}}$ for all
$i, j$, if  edge$_{ij}\ne 0$, else in $\sE_0$.

 We define
\eq{0.5}{H_*(X,D;E^\vee , \nabla^\vee) := H_*\Big(\sC_*(X;
E^\vee,\nabla^\vee)/\sC_*(D;E^\vee, \nabla^\vee)\Big) .}

(The growth condition means this depends on more than just the
topological sheaf $\sE^\vee$, so we keep $E^\vee,
\nabla^\vee$ in the notation.)

We now define a pairing
\begin{gather}\label{pairing}
(\ ,\ ):H^*_{DR}(X\setminus D; E,\nabla) \times H_*(X,D; E^\vee,
\nabla^\vee) \to 
\C;\quad *=0,1
\end{gather}
by integrating over chains in the following manner.
For $*=0$, then $H_0( X,D; E^\vee, \nabla^\vee)$ is generated by
 sections of the dual local system $\sE^\vee$ in points $\in X$
while $H^0_{DR}(X \setminus D; E,\nabla)$ is generated by global flat
sections in $\sE$ with moderate growth. So one can pair them.
For $*=1$, since $D \ne \emptyset$, then
$$H^1_{DR}(X\setminus D; E,\nabla)=H^0(X, \omega\otimes
E(*D))/\nabla H^0(X, E(*D)),$$ 
and since classes $c\otimes \epsilon$
generating $H_0( X,D; E^\vee, \nabla^\vee)$ have rapid decay, the integral
$\int_c <f_i c^*(e_i), \alpha>$ is convergent, where $\alpha 
\in H^0(X, \omega\otimes
E(*D))$ and $< \ >$ is the duality between $E^\vee$ and $E$.

The rest of the note is devoted to the proof of the following theorem.

\begin{thm}\label{mainthm} The process of integrating forms over
chains is compatible with homological and cohomological
equivalences and
defines a perfect pairing of finite dimensional vector complex spaces
$$(\ ,\ ):H^*_{DR}(X 
\setminus D; E,\nabla) \times H_*(X,D; E^\vee, \nabla^\vee) \to
\C;\quad *=0,1. 
$$
\end{thm}

\begin{ex}
\noindent {\rm (i).} If $\nabla$ has regular singular points, there are no
rapidly decaying flat sections, so $H_*(X,D; E^\vee, \nabla^\vee)\cong
H_*(X \setminus D; \sE^\vee)$. Also, $H^*_{DR}(X \setminus
D; E,\nabla)\cong H^*(U, \sE)$
(cf. \cite{DeI}, Th\'eor\`eme 6.2), 
and the theorem becomes the classical duality between
homology and cohomology. \newline
\rm{(ii).} Suppose $X=\P^1$, $D=\{0,\infty\}$. Let $E=\sO_{\P^1}$ with
connection $\nabla(1) = -dt + s\frac{dt}{t}$, for some $s\in
\C \setminus \{0,1,2,\ldots\}$. Then $\sE\subset E_U =
\sO_U$ is the trivial local system spanned by $e^tt^{-s}$, so
$\sE^\vee\subset E_U^\vee = \sO_U$ is spanned by $e^{-t}t^s$. We consider
the pairing $H^1_{DR} \times H_1 \to \C$ from theorem \ref{mainthm}. Note
first that $H^1_{DR}$ has dimension $1$, spanned by $\frac{dt}{t}$. This
can either be checked directly from \eqref{0.1}, using 
$$\nabla(t^p) =((p+s)t^{p-1}-t^p)dt,
$$ 
or by showing the de Rham cohomology is isomorphic
to the hypercohomology of the complex $\sO_{\P^1} \stackrel{\nabla}{\to}
\omega((0)+2(\infty))$, which is easily computed.  To compute $H_1(X,D;
E^\vee, \nabla^\vee)$, the singularity at $0$ is regular, so there are no
non-constant, rapidly decaying chains at $0$. The section
$\epsilon^\vee :=e^{-t}t^s$ of
$\sE^\vee$ is rapidly decaying on the positive real axis near $\infty$,
so the chain $c\otimes \epsilon^\vee$ in fig. $1$ above represents a
$1$-cycle. We have
$$(c\otimes e^{-t}t^s,\frac{dt}{t}) = (e^{2\pi is} - 1)\int_0^\infty
e^{-t}t^s \frac{dt}{t}
$$
which is a variant of Hankel's formula
(see \cite{W}, p. 245).
\newline
\noindent \rm{(iii).} Let $X, D, E$ be as in (ii), but take $\nabla(1)
=\frac{1}{2}(d(zu) - d(\frac{z}{u}))$ for some $z\in \C
\setminus \{0\}$.
 Here the connection has pole order $2$ at
$0$ and $\infty$ and it has trivial monodromy. Arguing as above, one
computes $\dim H^1_{DR} = 2$, generated by $u^pdu, \ p\in \Z$, with
relations $u^pdu = -\frac{2p}{z}u^{p-1}du - u^{p-2}du$. The Gau\ss-Manin
connection on this group is 
$$\nabla_{GM}(u^pdu) = \frac{1}{2}(u^{p+1}-u^{p-1})du\wedge dz.
$$
Assume $\text{Im}(z)>0$. Then the vector space
$H_1(\P^1,\{0,\infty\};E^\vee,\nabla^\vee)$ is generated by 
$$\{|u|=1\}\otimes \exp(\frac{1}{2}z(u-\frac{1}{u})), \quad
\text{and}\quad [0,i\infty]\otimes \exp(\frac{1}{2}z(u-\frac{1}{u})).
$$
(If $\text{Im}(z)\not > 0$, then the second path must be modified.) The
integrals
\begin{gather*}J_n(z) := \int_{\{|u|=1\}}
\exp(\frac{1}{2}z(u-\frac{1}{u}))\frac{du}{u^{n+1}}; \\
H_n(z):= \int_0^{i\infty}
\exp(\frac{1}{2}z(u-\frac{1}{u}))\frac{du}{u^{n+1}}
\end{gather*}
are periods and satisfy the Bessel differential equation
$$z^2\frac{d^2y}{dz^2} + z\frac{dy}{dz} + (z^2-n^2)y = 0
$$
The function  $J_n$ is entire.
To show that $H_n$ is linearly independent of $J_n$, it
will then be sufficient to show that $H_n$ is unbounded on
the positive part of the imaginary axis $\text{Re}(z)=0$ as $z\to 0$. 
Making the coordinate change $v=\frac{1}{u}$, and replacing $y$
by $\frac{1}{2}y$ one is led to show that $E_n(y)=\int_0^\infty
\exp(-y(v+\frac{1}{v}))\frac{dv}{v^{n+1}}$ is unbounded for
$y >0, y\to 0$. Writing $E_n(v)=\int_0^1+\int_1^\infty$, and
making the change of variable $v\to \frac{1}{v}$ in the integal
$\int_0^1$, one obtains 
\begin{gather*}
E_n(y)=\int_1^\infty \exp(-y(v+\frac{1}{v}))(\frac{1}{v^{n+1}}
+ v^{n-1})dv\\
\ge \int_1^\infty \exp(-2yv)(\frac{1}{v^{n+1}}
+ v^{n-1})dv.\\
\end{gather*}
For $|n| \ge 1$, then this expression is $\ge \int_1^\infty
\exp(-2yv)dv$ which is obviously unbounded. For $n=0$, one has
\begin{gather*}
E_0(y)\ge 2\int_1^\infty \exp(-2yv)\frac{dv}{v}\\
\ge 2\int_{2y}^\infty \exp(-v)\frac{dv}{v}\\
\ge 2\int_{2y}^1 \exp(-v)\frac{dv}{v},
\end{gather*}
where in the last inequality, we have assumed that $2y \le 1$. This last
integral is, upto something bounded, equal to $2\int_{2y}^1
\frac{dv}{v} =-2\log(2y)$, which is unbounded, as $y>0, y\to 0$.

Usually, for integers $n\in \Z$, 
one considers $J_n$ as one standard solution, but not
$H_n$ (see \cite{W}, p.371).
Finally, to
get Bessel functions for non-integral values of $n$, one may consider the
connection $\nabla(1)
=\frac{1}{2}(d(zu) - d(\frac{z}{u}))-n\frac{du}{u}$. 
\end{ex}

\section{Chains}

Let $D=\{x_1,\dotsc,x_n\}$ be as above, 
and let $\Delta_i$ be a small
disk about $x_i$ for each $i$. Let  $\delta_i$ be the boundary circle.
Define 
\ml{1.1}{H_*(\Delta_i,\delta_i\cup \{x_i\};E,\nabla) \\
=
H_*\Big(\sC_*(\Delta_i;E,\nabla)/(\sC_*(\delta_i;E,\nabla)+\sC_*(\{x_i\};
E,\nabla) \Big) }     (Note, for a set like $\delta_i$ which is closed and
disjoint from $D$, our chains coincide with the usual topological chains
with values in the local system $\sE$. The group 
$\sC_*(\{x_i\}; E,\nabla)$ consists of
constant chains $c:\Delta^n \to \{x_i\}$ with values in 
$$\sE_{x_i} := \sE_x/(1-\mu_i)\sE_x
$$
for some $x$ near $x_i$ as in \eqref{0.3}, where $\mu_i$ is the
local monodromy around $x_i$.) In the following theorem,
$H_*(U,\sE)$ is the standard homology associated to the local system
on $U=X\setminus D$. 
\begin{thm}\label{thm1.1} With notation as above, there is a long exact
sequence
\ml{1.2}{ 0 \to H_1(U, \sE) \to H_1(X,D; E, \nabla) \to \oplus_i
H_1(\Delta_i,\delta_i\cup \{x_i\};E,\nabla)  \\
\to H_0(U, \sE) \to H_0(X,D; E, \nabla) \to 0
}
\end{thm}    
\begin{proof}Let $\sC_*:= \sC_*(X; E, \nabla)/\sC_* (D;
E, \nabla)$  be the complex calculating
$H_*(X,D;E,\nabla)$, and let 
$$\sC_*(U)\subset \sC_*
$$
be the subcomplex calculating $H_*(U,\sE)$, i.e. the subcomplex of chains
whose support is disjoint from $D$. Of course, one has
$\sC_*(U;E,\nabla)=\sC_*(U; \sE)$, which justifies the notation.

Write $\sB = \sC_*/\sC_*(U)$. There
is an evident map of complexes 
\eq{1.3}{\psi:\oplus_i \sC_*(\Delta_i,\delta_i\cup 
\{x_i\};E,\nabla) \to \sB
}
which must be shown to be a quasi-isomorphism. 
Let
\begin{gather*}\sB(i)=\psi(\sC_*(\Delta_i,\delta_i\cup
\{x_i\};E,\nabla))=\\
\sC_*(\Delta_i,\delta_i\cup
\{x_i\};E,\nabla)/\sC_*(\Delta_i\setminus \{x_i\};
\sE)
 \subset
\sB.
\end{gather*}
Obviously the map $\alpha:\oplus_i \sB(i) \hookrightarrow
\sB$ is an inclusion. We claim first that $\alpha$ is a
quasi-isomorphism. To see this, note that all these complexes admit
subdivision maps $\text{subd}$ which are homotopic to the identity. Given
a chain $c\in \sB$, there exists an $N$ such that $\text{subd}^N(c) \in
\oplus \sB(i)$. Taking $c$ with $\p c = 0$, it follows that
$\oplus H_*(\sB(i))$ surjects onto $H_*(\sB)$. If $\alpha(x) = \p y$, we
choose $N$ such that $\text{subd}^N(y) = \alpha(z)$. Since $\alpha$ is
injective and commutes with $\text{subd}$, it follows that $\alpha$ is
injective on homology as well, so $\alpha$ is a quasi-isomorphism. 

It remains to show the surjective map of complexes 
$$\beta:
\sC_*(\Delta_i,\delta_i\cup \{x_i\};E,\nabla) \to \sB(i)
$$
is a quasi-isomorphism. The kernel of $\beta$ is
$$\sC_*(\Delta_i\setminus \{x_i\}; \sE)/\sC_*(\delta_i; \sE),$$
which is acyclic as
$\delta_i \hookrightarrow \Delta_i \setminus \{x_i\}$ 
admits an evident homotopy
retract. 

The next point is to show
\eq{1.4}{H_*(\Delta_i,\delta_i\cup \{x_i\};E,\nabla)=(0); \ i=0,2.
}
The assertion for $H_0$ is easy because any point $y$ in $\Delta_i \setminus
\{x_i\}$ can be attached to $\delta_i$ by a radial path $r$ not passing
through $x_i$. Then $\epsilon \in \sE_y$ extends uniquely to $\epsilon$
on $r$ and $\p(r\otimes \epsilon) \equiv y\otimes \epsilon \mod
\text{chains on $\delta_i$}$. Vanishing in \eqref{1.4} when $i=2$ will be
proved in a sequence of lemmas. For convenience we drop the subscript
$i$ and replace $x_i$ with $0$. 
\begin{lem}\label{lem1.2} Let $\ell\subset \Delta$ be a radial line
meeting
$\delta$ at
$p$. Let $\sE_\ell$ be the space of sections of the local system along
$\ell\setminus \{0\}$ with rapid decay at $0$. Then
$$H_*(\ell,\{0,p\};E,\nabla) \cong \begin{cases}0 & *\neq 0 \\
\sE_\ell & *=1\end{cases}
$$
\end{lem}
\begin{proof}[proof of lemma] Let $\sC_*(\ell)$ be the complex of chains
calculating this homology, and let $\sC_*(\ell\setminus \{0\})
\subset \sC_*(\ell)$ be 
the subcomplex of chains not meeting $0$. Then $\sC_*(\ell\setminus\{0\})$ is
contractible, and 
$$\sC_*(\ell)/\sC_*(\ell\setminus\{0\}) \cong (\bold C_*(\ell)/\bold
C_*(\ell\setminus \{0\}))\otimes \sE_\ell
$$
where $\bold C_*$ denotes classical topological chains. The result follows.
\end{proof}

One knows from the theory of irregular connections in $\dim 1$ \cite{M1}
that $\Delta \setminus \{0\}$ 
can be covered by open sectors $V\subsetneq\Delta$
such than 
\eq{1.5}{E,\nabla |_V \cong \oplus_i (L_i\otimes M_i)
}
where $L_i$ is rank $1$ and $M_i$ has a regular singular point. Let $W
\subset V\cup \{0\}$ be a smaller closed sector with outer boundary
$\delta_W = \delta\cap W$ and radial sides $\ell_1,\ell_2$. Recall the
Stokes lines are radial lines where the horizontal sections of the $L_i$
shift from rapid decay to rapid growth. We assume $W$ contains at most
one Stokes line, and that $\ell_1, \ell_2$ are not Stokes lines. 
Writing $W=W_1\cup W_2$, where $W_i$ are even smaller sectors,
each of which containing the Stokes line if there is one, one
may think of the following lemma as a Mayer-Vietoris sequence.
\begin{lem}\label{lem1.3} With notation as above, Let $w$ be a basepoint
in the interior of $W$. then 
$$H_*(W,\delta_W\cup \{0\};E,\nabla) \cong \begin{cases} 0 & *\neq 1 \\
\sE_{\ell_1} + \sE_{\ell_2} \subset \sE_w & *=1
\end{cases}
$$
\end{lem}
\begin{proof}[proof of lemma] One has 
$$\oplus _iH_1(\ell_i,\{0,p_i\};E,\nabla) \to H_1(W,\delta_W\cup
\{0\};E,\nabla)
$$
and of course the assertion of the lemma is that this coincides with
$\sE_{\ell_1}\oplus \sE_{\ell_2} \to \sE_{\ell_1}+ \sE_{\ell_2}$. To
check this, by \eqref{1.5} one is reduced to the case $E=L\otimes M$
where $L$ has rank $1$ and $M$ has regular singular points. 

If $W$ does not contain a Stokes line for $L$ then $\sE_{\ell_1} =
\sE_{\ell_2}= \sE_{\ell_1}+\sE_{\ell_2}$, and the argument is exactly as
in lemma \ref{lem1.2}. 

Suppose $W$ contains a Stokes line for $L$. Then (say) $\sE_{\ell_1} =
\sE_w$ and  $\sE_{\ell_2} = (0)$. Let $\sC_*(W)$ be the complex of chains
calculating the desired homology, and let $\sC_*(W\setminus
\{0\})\subset \sC_*(W)$ be
the chains not meeting $0$. As in the previous lemma,
$\sC_*(W\setminus \{0\})$ is
acyclic. We claim the map
$$\sC_*(\ell_1) \to \sC_*(W)/\sC_*(W\setminus\{0\})
$$
is a quasi-isomorphism. If we
choose an angular coordinate $\theta$ such that 
$$\ell_1 : \theta = 0;\quad \text{Stokes} : \theta = a>0;\quad \ell_2 :
\theta = b>a,
$$
then rotation $re^{i\theta} \mapsto re^{(1-t)i\theta}$ provides a
homotopy contraction of the inclusion of $\ell_1 \subset W$. This
homotopy contraction preserves the condition of rapid decay, proving the
lemma.
\end{proof}

Let $\pi_d:\Delta \to \Delta$ be the ramified cover of degree $d$
obtained by taking the $d$-th root of a parameter at $0$.  By the theory
of formal connections
\cite{M1}, one has, for suitable $d$, a decomposition as in \eqref{1.5}
for the formal completion of the pullback
$\widehat{\pi_d^*E}\cong \oplus_i L_i\otimes M_i$. Let $m_i$ be the degree
of the pole of the connection on
$L_i$ when we identify $L_i \cong \widehat{\sO}$, i.e. $\nabla_{L_i}(1) =
g_i(z)dz$ for a local parameter $z$, and $m_i$ is the order of pole of
$g_i$. 
\begin{lem}\label{lem1.4} We have
$$\dim H_p(\Delta,\delta\cup \{0\}; E,\nabla) = \begin{cases} 0
& p\neq 1 \\
\frac{1}{d}\sum_{m_i\ge 2} (m_i-1)\dim(M_i) & p=1.\end{cases}
$$
\end{lem}
\begin{proof}[proof of lemma]Assume first that we have a decomposition of
the type \eqref{1.5} on $\widehat{E}$ itself, i.e. that no pullback
$\pi_d^*$ is necessary. We write
$\Delta$ as a union of closed sectors $W_0,\dotsc,W_{N-1}$ where $W_i$
has radial boundary lines
$\ell_i$ and $\ell_{i+1}$. We assume each $W_i$ has at most one Stokes
line. Using excision together with the previous lemmas we get
\ml{1.6}{0 \to H_2(\Delta,\delta\cup \{0\}; E,\nabla) \to
\oplus_{i=0}^{N-1} H_1(\ell_i, \{p_i,0\};E,\nabla)  \\
\stackrel{\nu}{\to}\oplus_{i=0}^{N-1}H_1(W_i,\delta_{W_i}\cup \{0\};
E,\nabla)
\to H_1(\Delta,\delta\cup \{0\}; E,\nabla) \to 0
 }
By lemma \ref{lem1.3}, the map $\nu$ above is given by
$$ \nu(e_0,\dotsc,e_{N-1}) = (e_0-e_1,e_1-e_2,\dotsc,e_{N-1}-e_0).
$$
An element in the kernel of $\nu$ is thus a section $e$ of
$\sE|_{\Delta-\{0\}}$ which has rapid decay along each $\ell_i$. Since each
$W_i$ contains at most one Stokes line, such an $e$ would necessarily
have rapid decay on every sector and thus would be trivial. This proves
vanishing for $H_2(\Delta,\delta;E,\nabla)$. Finally, to compute the
dimension of $H_1$, note that if $L_i$ has a connection with pole of
order $m_i$, then it has a horizontal section of the form $e^f$, where
$f$ has a pole of order $m_i-1$. (The connection is $1\mapsto df$.)
Suppose $f=az^{1-m_i} +\ldots$. Stokes lines for this factor are radial
lines where $az^{1-m_i}$ is pure imaginary. Thus, there are $2(m_i-1)$
Stokes lines for this factor. Consider one of the Stokes lines, and
suppose it lies in $W_k$. If the real part of $az^{1-m_i}$ changes from
negative to positive as we rotate clockwise through this line, say we are
in case $+$, otherwise we are in case $-$. We have
\eq{1.7}{\dim(\sE_{\ell_k} + \sE_{\ell_{k+1}}) - \dim \sE_{\ell_k} =
\begin{cases} 0 & \text{case } + \\
\dim(M_i) & \text{case }-, \end{cases} 
}
since the two cases alternate, we get a contribution of
$(m_i-1)\dim(M_i)$. If $m_i \le 1$ there are no rapidly decaying
sections, so that case can be ignored. Summing over $i$ with $m_i\ge 2$
gives the desired result. 

Finally, we must consider the general case when the decomposition
\eqref{1.5} is only available on $\widehat{\pi_d^*E}$ for some $d\ge 2$.
By a trace argument, vanishing of the homology upstairs, i.e. for
$\widehat{\pi_d^*E}$, in degrees $\neq 1$ implies vanishing downstairs. Since
$\pi_d: \Delta\setminus \{0\} \to \Delta \setminus\{0\}$ 
is unramified, an Euler
characteristic argument (or, more concretely, just cutting into small
sectors over which the covering splits) shows that the Euler
characteristic multiplies by $d$ under pullback, proving the lemma. 
\end{proof}

In particular, we have now completed the proof of theorem \ref{thm1.1}. 
\end{proof}

\section{de Rham Cohomology}

In this section, using differential forms, we construct the dual sequence
to the homology sequence from theorem \ref{thm1.1}. (More precisely, we
continue to work with $E,\nabla$, so the sequence we construct will be
dual to the homology sequence with coefficients in $E^\vee,
\nabla^\vee$). Consider the diagram of complexes
\eq{2.1}{\minCDarrowwidth.5cm\begin{CD}0 @>>> E(*D) @>>> j_*E_U  @>>>
j_*E_U/E(*D) @>>> 0 \\
@. @V\nabla_{{\rm mero}} VV @V\nabla_{{\rm an}} VV @V\nabla_{{\rm an/mero}} VV \\
0 @>>> E(*D)\otimes \omega @>>> j_*E_U\otimes \omega  @>>>
(j_*E_U/E(*D))\otimes \omega @>>> 0
\end{CD} 
}
A result of Malgrange \cite{M3} is that $\nabla_{{\rm an/mero}}$ is surjective.
Define $N:= \oplus_i N_i = \ker(\nabla_{{\rm an/mero}})$. Since none of these
sheaves has higher cohomology (by assumption $D\neq \emptyset$) we get a
5-term exact sequence by taking global sections and applying the serpent
lemma:
\ml{2.2}{0 \to H^0_{DR}(U;E,\nabla)  \to H^0(U,\sE) \to N \\
\to H^1_{DR}(U;E,\nabla)  \to H^1(U,\sE) \to 0
}

\begin{thm}\label{thm2.1} Integration of forms over chains defines a
perfect pairing between the exact sequence \eqref{2.2} and the exact
sequence from theorem \ref{thm1.1}:
\ml{2.3}{ 0 \to H_1(U, \sE^\vee) \to H_1(X,D; E^\vee, \nabla^\vee) \to
\oplus_i H_1(\Delta_i,\delta_i\cup \{x_i\};E^\vee,\nabla^\vee)  \\
\to H_0(U, \sE^\vee) \to H_0(X,D; E^\vee, \nabla^\vee) \to 0.
}
\end{thm}
\begin{proof}To establish the existence of a pairing, note that if
$c\otimes \epsilon^\vee$ is a rapidly decaying chain and $\eta$ is a form
of the same degree with moderate growth, then elementary estimates show
$\int_c \langle\epsilon^\vee,\eta\rangle$ is well defined. Suppose
$c:\Delta^n
\to X$ and write $\Delta^n = \lim_{t\to 0}\Delta^n_t$ where $\Delta^n_t$
denotes
$\Delta^n \setminus 
\text{tubular neighborhood of radius $t$ around $\p \Delta^n$}$.
 Let $c_t = c|_{\Delta^n_t}$ and suppose $\eta = d\tau$ where $\tau $ has
moderate growth also. Then
\eq{2.4}{ \int_c\langle\epsilon^\vee,\eta\rangle = \lim_{t\to 0}
\int_{c_t} \langle \epsilon^\vee, d\tau \rangle = \lim_{t\to 0}
\int_{\p c_t} \langle \epsilon^\vee,  \tau \rangle = \int_{\p c}\langle
\epsilon^\vee,  \tau \rangle.   
}
Note $\p c$ may include simplices mapping to $D$. Our definition
\eqref{0.5} of $\sC_*(X,D;E^\vee,\nabla^\vee)$
factors these chains out. Thus, we do get a pairing of complexes. 

Of course, chains away from $D$ integrate with forms with possible
essential singularities on $D$. To complete the description of the
pairing, we must indicate a pairing
\eq{}{(\ ,\ ):N_i \times H_1(\Delta_i,\delta_i\cup
\{x_i\};E^\vee,\nabla^\vee) 
\to \C. 
}
To simplify notation we will drop the subscript $i$ and take $x_i=0$.  An
element in $H_1$ can be represented in the form $\epsilon^\vee\otimes c$
where $c$ is a radial path. Let $c\cap \delta = \{p\}$. Given $n\in N$,
choose a sector $W$ containing $c$ on which $\sE$ has a basis
$\epsilon_i$. By assumption, we can represent $n=\sum a_i\epsilon_i$ with
$a_i$ analytic on the open sector, such that
\eq{2.6}{\nabla(\sum a_i\epsilon_i) = \sum \epsilon_i\otimes da_i = \sum
e_i\otimes \eta_i 
}
where $e_i$ from a basis of $E$ in a neighborhood of $0$ and $\eta_i$ are
meromorphic $1$-forms at $0$. then by definition
\eq{2.7}{(\epsilon^\vee \otimes c,n) := \int_c \sum_i \langle
\epsilon^\vee,e_i \rangle \eta_i - \sum_i \langle
\epsilon^\vee,\epsilon_i \rangle a_i(p).  
}
The pairing is taken to be trivial on chains which do not contain $0$. If
$s$ is a $2$-chain bounding two radial segments $c$ and $c'$ and a path
along $\delta$ from $p$ to $p'$. Then Cauchy's theorem (together with a
limiting argument at $0$) gives
\ml{2.8}{0 = \int_c \sum_i \langle
\epsilon^\vee,e_i \rangle \eta_i - \int_{c'} \sum_i \langle
\epsilon^\vee,e_i \rangle \eta_i + \int_p^{p'}\sum_i \langle
\epsilon^\vee,\epsilon_i \rangle da_i \\
= (\epsilon^\vee \otimes c,n) - (\epsilon^\vee \otimes c',n).
}
Similar arguments show the pairing independent of the choice of the
radius of the disk. Also, if $\sum a_i\epsilon_i = \sum b_i e_i$ with
$b_i$ meromorphic at $0$, then 
\ml{2.9}{(\epsilon^\vee \otimes c,n) = \int_c \sum \langle
\epsilon^\vee,e_i \rangle \eta_i - \sum_i \langle \epsilon^\vee,
\epsilon_i \rangle a_i(p) \\
= \int_c d\langle\epsilon^\vee,\sum b_i e_i \rangle - \sum_i \langle 
\epsilon^\vee, \epsilon_i \rangle a_i(p)   \\
= \langle \epsilon^\vee, \sum b_i e_i \rangle(p) - \sum_i \langle 
\epsilon^\vee, \epsilon_i \rangle a_i(p) = 0
}
It follows that the pairing is well defined.

\begin{lem}The diagrams 
$$\begin{array}{ccc}H_1(X,D;E^\vee, \nabla^\vee) & \to & \oplus
H_1(\Delta_i,\delta_i;E^\vee, \nabla^\vee) \\
\times && \makebox[3.5cm][l]{$\times$} \\
H^1(X\setminus D;E,\nabla) & \leftarrow & \makebox[3.5cm][l]{$\oplus N_i$} \\
\makebox[3.8cm][r]{$\searrow$} && \makebox[4cm][l]{$\swarrow$}
\rule{0cm}{.7cm} \\
& \makebox[1cm][c]{$\C$}
\end{array}
$$
and 
$$\begin{array}{ccc} \oplus
H_1(\Delta_i,\delta_i;E^\vee, \nabla^\vee) & \to & H_0(U, \sE^\vee) \\
\makebox[3cm][r]{$\times$} && \times \\
\makebox[3cm][r]{$\oplus N_i$} & \leftarrow & H^0(U, \sE) \\
\makebox[3.5cm][r]{$\searrow$} &&\swarrow \qquad \rule{0cm}{.7cm} \\
& \C 
\end{array}
$$
commute. 
\end{lem}
\begin{proof}[proof of lemma] Consider the top square. The top arrow is
excision, replacing a chain with the part of it lying in the disks
$\Delta_i$. The bottom arrow maps an $n$ as above in some $N_i$ to $\sum
e_j\otimes
\eta_j = \sum \epsilon_j\otimes da_j$. Along $c$ outside the disks $\sum
e_j\otimes \eta_j$ is exact; its integral along the chain is a sum of
terms of the form  $\sum_i \langle 
\epsilon^\vee, \epsilon_j \rangle a_j(p_i)$ where $p_i\in c\cap
\delta_i$. For the part of the chain inside the $\Delta_i$ of course we
must take $\int_{c\cap \Delta_i}\langle \epsilon^\vee,e_j\rangle\eta_j$.
Combining these terms with appropriate signs yields the desired
compatibility. 

For the bottom square, the top arrow associates to a relative chain on
$\Delta_i$ its boundary on $\delta_i \subset U$. The bottom arrow
associates to a horizontal section $\epsilon$ on $U$ the corresponding
element in $N$. Note here the $a_j$ will be constant so in the pairing
with $N$ only the term $-\sum \langle
\epsilon^\vee,\epsilon_j\rangle a_j(p)$ survives. The assertion of the
lemma follows. 
\end{proof}

Returning to the proof of the theorem, we see it reduces to a purely
local statement for a connection on a disk. In the following lemma, we
modify notation, writing $N$ to denote the corresponding group for a
connection on a disk $\Delta$ with a meromorphic singularity at $0$. 
\begin{lem} The pairing
$$(\ ,\ ):N \times H_1(\Delta,\delta;E^\vee, \nabla^\vee) \to \C
$$
is nondegenerate on the left, i.e. $(\epsilon^\vee \otimes c,n) = 0$ for
all relative $1$-cycles implies $n=0$. 
\end{lem}
\begin{proof}[proof of lemma] We work in a sector and we suppose the basis
$\epsilon_i$ taken in the usual way compatible (in the sector) with the
decomposition into a direct sum of rank $1$ irregular connections tensor
regular singular point connections. Let $\epsilon_i^\vee$ be the dual
basis. 

Fix an $i$ and suppose first $\epsilon_i$ and $\epsilon_i^\vee$ both have
moderate growth. We claim $a_i$ has moderate growth. For
this it suffices to show $da_i$ has moderate growth. But
\eq{2.10}{ da_i = \langle \nabla(n),\epsilon_i^\vee \rangle =
\sum_j\langle e_j,\epsilon_i^\vee \rangle \eta_j. 
}
This has moderate growth because, $e_j, \epsilon_i^\vee,\text{and }
\eta_j$ all do. 

Now assume $(\epsilon^\vee \otimes c,n) = 0$ for all $\epsilon^\vee
\otimes c \in H_1$. Fix an $i$ and assume $\epsilon_i^\vee$ is rapidly
decreasing in our sector. Let $c$ be a radius in the sector with
endpoint $p$. We can find (cf. \cite{M1}, chap. IV, p.53-56) a basis
$t_i$ of
$E$ on the sector with moderate growth and such that $ t_i=
\psi_i\epsilon_i$, so $t_i^\vee = \psi_i^{-1}\epsilon_i^\vee$. 

We are interested in the growth of $a_i\epsilon_i$ along $c$. We have
\eq{2.11}{a_i(p)\epsilon_i(p) = \Big(\int_c\sum_j\langle
\epsilon_i^\vee,e_j \rangle \eta_j\Big)\epsilon_i(p) =  \Big(\int_c
\psi_i\sum_j\langle t_i^\vee,e_j \rangle
\eta_j\Big)\psi_i(p)^{-1}t_i(p). 
}
Asymptotically, taking $y$ the parameter along $c$, $\psi_i(y)\sim
\exp(-ky^{-N})$ as $y\to 0$ for some $k>0$ and some $N\ge 1$. We need to
know the integral
\eq{2.12}{ \exp(kp^{-N})\int_0^p y^{-M}\exp(-ky^{-N})dy
}
has moderate growth as $p \to 0$. Changing variables, so $x=y^{-1},
q=p^{-1}, u=x-q$, this becomes
\ml{2.13}{ \int_0^\infty (u+q)^{M-2} \exp(q^N - (u+q)^N)du \\
= \int_0^\infty
(u+q)^{M-2} \exp(-u^N - qf(u,q))du, 
}
where $f$ is a sum of monomials in $q$ and $u$ with positive
coefficients.  Clearly this has at worst polynomial growth as $q\to
\infty$ as desired. 

Finally, assume $\epsilon_i^\vee$ is rapidly increasing and
$\epsilon_i$ is rapidly decreasing. We have as above
\eq{2.14}{\sum_i e_j\otimes \eta_j = \sum_j \epsilon_j\otimes da_j =
\sum_j
\psi_j^{-1}t_j\otimes da_j 
}
In particular, $\psi_i^{-1}da_i$ has moderate growth. This implies
$a_i\epsilon_i=a_i\psi_i^{-1}t_i$ has moderate growth as well. Indeed,
changing notation, this amounts to the assertion that if $g$ is rapidly
decreasing and $g\frac{df}{dz}$ has moderate growth, then $gf$ has
moderate growth. Fix a point $p_0$ with $0<p<p_0$. the mean value theorem
says there exists an $r$ with $p\le r\le p_0$ such that
$$g(p)f(p) = g(p)(f(p_0) + (p-p_0)f'(r))
$$
Suppose $|f'(q)g(q)|<<q^{-N}$. We get
$$|g(p)f(p)| << |g(p)f'(r)| \le|g(r)f'(r)|<< r^{-N}\le p^{-N}
$$
proving moderate growth. 

We conclude that our representation for $n$ has moderate growth, and
hence it is zero in $N$. It follows that the pairing $N\times H_1 \to \C$
is nondegenerate on the left. 
\end{proof}

Returning to the global situation, we have now 
$$\dim N_i \le \dim
H_1(\Delta_i, \delta_i;E^\vee, \nabla^\vee), 
$$
and to finish the proof of the theorem, it will suffice to show these
dimensions are equal. 

\begin{lem} With notation as above, $\dim N_i = \dim H_1(\Delta_i,
\delta_i;E^\vee, \nabla^\vee)$. 
\end{lem}
\begin{proof}[proof of lemma] It will suffice to compute the difference
of the two Euler characteristics
\eq{2.15}{ \chi(U, \sE) - \chi_{DR}(U;E,\nabla). 
}
It is straightforward to show this difference is invariant if $U$ is
replaced by a smaller Zariski open set, and that the Euler
characteristics are multiplied by the degree in a finite \'etale covering
$V \to U$. Using lemma \ref{lem1.4}, we reduce to the case where
formally locally at each $x_i\in D$ we have $E\otimes \widehat{K}_{x_i}
\cong \oplus_j L_{ij}\otimes M_{ij}$ with $L_{ij}$ rank $1$ and $M_{ij}$
at worst regular singular. (Here $\widehat{K}_{x_i}$ is the Laurent power
series field at $x_i$). Let $m_{ij}$ be the degree of the pole for the
connection on $L_{ij}$. Then one can find coherent sheaves 
$$F_2 \subset F_1 \subset E(*D)
$$ 
such that 
\begin{gather*}F_1/F_2 \cong \oplus_{ij}M_{ij}/M_{ij}(-m_{ij}x_i) \\
E^\nabla \subset H^0(F_2); \nabla(F_2) \subset F_1\otimes \omega \\
H^0(F_1\otimes \omega) \twoheadrightarrow H^1_{DR}(U;E,\nabla) 
\end{gather*}
It follows that, writing $g=\text{genus}(X)$
\ml{2.16}{\chi_{DR}(U;E,\nabla) \\ 
= \chi(F_2) - \chi(F_1\otimes \omega) =
-\text{rk}(E)(2g-2) - \sum_{ij} m_{ij}\dim(M_{ij}).
}
Since
\eq{2.17}{\chi(U,\sE) = -\text{rk}(E)(2g-2+n), 
}
(which is proven algebraically as above, replacing $\nabla$ by
the regular connection associated to $\sE$)
it follows that
$$\chi_{DR}(U;E,\nabla) - \chi(U,\sE) =- \sum_{ij} (m_{ij}-1)\dim(M_{ij})
$$
Referring to lemma \ref{lem1.4}, 
we see that this is the desired formula. 
\end{proof}

This completes the proof of the theorem. 
\end{proof}

\newpage
\bibliographystyle{plain}
\renewcommand\refname{References}

\end{document}